\documentclass[12pt]{article}
\usepackage[russian]{babel}
\usepackage[cp1251]{inputenc}
\usepackage{amsmath,amssymb}
\begin{document}
\thispagestyle{empty}
\begin{center}

{\bf  ЭРГОДИЧЕСКИЕ ВОЛЬТЕРРОВСКИЕ КВАДРАТИЧНЫЕ ОТОБРАЖЕНИЯ СИМПЛЕКСА }\\
\vspace{0.5cm}

{\bf Н.Н.Ганиходжаев, Д.В.Занин}\\[0.4mm]
\end{center}

 {\bf Аннотация.}
 В данной работе дается полный анализ асимптотического поведе\-ния
 траекторий вольтерровских квадратичных отображений трехмерного симплекса в себя.

{\bf Annotation}
In the paper a Volterra quadratic stochastic
operators of three dimensional simplex into itself is
considered.The full description of ergodic properties such
operators is given.

{\bf Ключевые слова }:Квадратичный стохастический оператор;
симплекс; эргодическая теорема.

\section{Введение}
В математической биологии (теории наследования) и в теории
вероятностей возникает вопрос об изучении траекторий одного класса
нелинейных преобразо\-ваний называемых квадратичными
стохастическими операторами [1-9,11-14].

Квадратичный стохастический
оператор (ксо), отображающий стандартный симплекс

$$ S^{m-1}=\{\textbf{x}=(x_1,...,x_m)\in R^m: x_i\geq 0, \sum^m_{i=1}x_i=1 \} \eqno(1)$$
в себя, имеет вид
$$ V: (V\textbf{x})_k=\sum^m_{i,j=1}p_{ij,k}x_ix_j, \ \ (k=1,...,m), \eqno(2)$$
где  коэффициенты $p_{ij,k}-$ удовлетворяют следующим условиям:
$$ p_{ij,k}\geq 0, p_{ij,k}=p_{ji,k}, \ \ \sum^m_{k=1}p_{ij,k}=1, \ \ (i,j,k=1,...,m). \eqno(3)$$
Заметим, что условие симметричности не является обременительным ,так
как в противном случае оно может быть получено переопределением
$p_{ij,k}$ следующим образом
$$ q_{ij,k}=\frac{p_{ij,k}+p_{ji,k}}{2}.$$
Легко проверить,что оператор не меняется при таком
переопределении.
 Следующая модель поясняет биологический смысл
ксо. Рассмотрим некоторую биологическую популяцию ,т.е. замкнутое
относительно размножения сообщество организмов.

Предположим, что каждая особь входящая в популяцию, принадлежит
некоторой единственной из $m $ разновидностей, $m =2,3,...$ Шкала
разновидностей должна быть такой,чтобы разновидности родителей
$(i,j)$ однозначно определяли вероятность каждой разновидности $k$
для непосредственного потомка первого поколения.Обозна\-чим эту
вероятность (коэффициент наследственности)через $p_{ij,k}$.Очевидно
в этом случае выполнены условия $ p_{ij,k}\geq 0, \ \
\sum^m_{k=1}p_{ij,k}=1, \ \ (i,j,k=1,...,m). $ Предполо\-жим
 что популяция настолько велика, что можно пренебречь флюкту\-ациями
частот.Тогда ее состояние можно описывать вектором
$\textbf{x}=(x_1,...,x_m)\in S^{m-1}$,где $x_i$ - доля
разновидности $i$ в популяции.В этом случае по формуле полной
вероятнос\-ти, распределение разновидностей в следующем поколении
задается вектором $V\textbf{x}\in S^{m-1}$ (2).
  Одной из основных проблем является задача изучения эргодических свойств таких операторов.
Рассмотрим классическое определение эргодичности топологической
динамической системы.

{\bf Определение 1.} Оператор $V$ называется эргодическим если и
только если для всякой непрерывной функции определенной на симплексе
и для любой начальной точки обшего положения $ \textbf{x}\in
S^{m-1}$ существует предел
   $$ \lim _{k\rightarrow \infty} \frac{1}{k}\sum_{n=0}^{k-1}  f(V^n \textbf{x}).$$

   Здесь понятие точки общего положения не конкретизируется, но имеется в
виду либо принадлежность к некоему множеству полной лебеговой
меры, либо к некоему множеству второй бэровской категории. Это
понятие введено Уламом [14]. Он рассматривал топологический аналог
теоретико- мерной эргодичности. В [13] на основе численных
расчетов была сформулирована гипотеза , что любой к.с.о.эргодичен.
Д.Занин и М.Юлдашев доказали(в неопубликованной дипломной работе)
справедливость гипотезы Улама при $ m=2$.В 1977 году Захаревич
доказал, что в общем случае эта гипотеза неверна.Им были построены
примеры неэргодических к.с.о. при $m=3$ [9].
 Ниже будем изучать поведение вольтерров\-ских стохастических операторов (названых в
честь V. Volterra, который рассмотрел эволюционные уравнения
биологии).

 {\bf Определение 2.} К.с.о. называется вольтерровским
если $ p_{ij,k}=0$ для всех $k\notin \{i,j\}$.

Биологический смысл этого определения достаточно прост- потомок
повторяет одного из родителей.

Вольтерровский оператор можно записать в виде
$$(V\textbf{x})_k=x_k(1+\sum_{i=1}^m a_{ki}x_i),$$
где $A=(a_{ij})_{i,j=1}^m $ -кососимметрическая матрица с
$a_{ij}=2p_{ji,j} -1$ и $|a_{ij}| \leq 1$. При $m=3$ произвольный
вольтерровский ксо может быть представлен следующим образом
$$ V:(x,y,z)\rightarrow (x(1+ax-bz),y(1-ax+cz),z(1+bx-cy)). \eqno (4)$$
В [3] установлено необходимое и
достаточное условие эргодичности вольтерровских ксо определенных
на $S^2$.

{\bf Теорема 1.}Вольтерровский ксо (4) неэргодичен тогда и только
тогда, когда параметры $a,b,c$ отличны от нуля и имеют один и тот же
знак.

Заметим ,что при $a=b=c=1$  получим пример рассмотренный
Захаревичем [9], а при $a=b=c=\varepsilon, 0<\varepsilon\leq 1$
получим семейство операторов указанных в Замечании 2 [9].

 В  данной работе рассмотрим случай размерности $m=4$.

\section{Сужение аттрактора}

Для произвольного вольтерровского оператора определенного на $
S^{4-1}$ ,соответст\-вующую кососимметрическую матрицу
$A=(\tilde{a}_{ij})_{i,j=1}^4 $
 запишем в следующем виде
 \begin{displaymath}
 \left(
 \begin{array} {cccc}
 0 & a_{12}& a_{13} & -a_{14}\\
 -a_{12} & 0 & a_{23} & a_{24} \\
 -a_{13} & -a_{23} & 0 & a_{34} \\
 a_{14} & -a_{24} & -a_{34} & 0
 \end {array}\right )\eqno (5)
 \end{displaymath}
Разобьем множество таких матриц на следующие три класса.К первому
классу отнесем матрицы содержащие строку, все элементы  которой
неотрицательны. Ко второму классу отнесем матрицы не содержащие
строку, все элементы которой неотрицательны, но содержащую строку,
все элементы которой неположительны.

Матрицы не попавшие в эти два класса, отнесем к третьему классу.
Легко видеть, что для матриц из третьего класса (изменяя, если надо,
порядок на множестве {1, 2, 3, 4}) все числа $ a_{ij}$ (а не
элементы матрицы) можно считать неотрицательными.
 Рассмотрим эти случаи по порядку.

\subsection{Матрицы первого типа}

Пусть все элементы $i$-й строки матрицы (5) неотрицательны. Тогда
для любого $ \textbf{x}=(x_1,x_2,x_3,x_4)\in S^{4-1} $ справедливо
следующее неравенство:
 $$(V \textbf{x})_i\geq x_i .$$
 Итак, последовательность $(V^k \textbf{\textbf{x}})_i$ не убывает. Пусть $\alpha$
предел этой последователь\-ности. Если $x_i$ > 0, то $\alpha$ > 0.
Перейдем к пределу при $k\rightarrow\infty$ в равенстве
$$(V^{k+1}\textbf{x})_i = (V^{ k}\textbf{x})_i(1 +\sum_{l\ne i}a_{il}(V^ {k}\textbf{x})_l).$$
 Имеем
 $$\lim_{k\rightarrow\infty}(V^ {k+1}\textbf{x})_i
= \lim _{k\rightarrow\infty}(V^{k}\textbf{x})_i(1 +\sum
a_{il}(V^{k}\textbf{x})_l).$$ Так как левая часть стремится к
$\alpha$ и множитель в правой части также стремится к $\alpha$, то
существует предел
$$\lim_{k\rightarrow\infty} \sum_{l\ne i}a_{il}(V^{k}\textbf{x})_l = 0.$$
 Пусть некоторые из $a_{il}$ строго положительны.
Тогда соответствующая координата вектора $(V^{k}\textbf{x})_l$
стремится к нулю. Предположим для простоты, что все они (кроме $l =
i$) строго положительны. Тогда, $V^{k}\textbf{x}$ стремится к
вершине симплекса. Отсюда следует регулярность а значит и
эргодичность оператора $V$.

\subsection{Матрицы второго типа}

 Пусть все элементы $i$-й строки матрицы (5) неположительны.Тогда для любого $ \textbf{x}=(x_1,x_2,x_3,x_4)\in
 S^{4-1} $
 $$ (V \textbf{\textbf{x}})_i\leq x_i$$
так что последовательность $(V^{k}\textbf{x})_i$ не возрастает и
поэтому у нее есть предел. Обозначим этот предел через $\beta$.
Предположим, что $\beta$> 0. Как и в предыдущем случае, получаем
существование предела
$$\lim_{k\rightarrow\infty} \sum_{l\ne i}a_{il}(V^{k}\textbf{x})_l = 0.$$
Предположим для простоты , что все $a_{il} < 0$ (кроме $l = i$).
Тогда траектория сходится к  $i$-ой вершине симплекса. Но $$\lim_{
k\rightarrow\infty} (V^{k}\textbf{x})_i\leq x_i$$  и $$\lim_{
k\rightarrow\infty}(V^{k}\textbf{x})_i = 1.$$ Это означает, что
начальная точка также была вершиной. Таким образом, наше
предположение что $\beta> 0$ приводит к тривиальному случаю.Пусть
$\beta =0$,то есть
$$ \lim_{k\rightarrow\infty}(V^{k}\textbf{x})_i = 0.$$
Без потери общности, полагаем $i = 4.$ Тогда симплекс
$\Gamma_{123}$ содержит внутри седловую точку, которую обозначим
через $ M_0$. Пусть $\omega$ предельное множество начальной точки
М.Очевидно оно содержится в $\Gamma_{123}$.Рассмотрим следующиe
два возможных случая.

Первый случай:  $\omega$-предельное множество точки M содержит
точку $ M_0$. Это означает, что для любой окрестности U точки $
M_0$ существует подпоследователь\-ность $n_k$, стремящаяся к
бесконечности, такая что $V^{n_k}\textbf{x}\in U.$ Определим
канони\-ческие координаты в окрестности $ M_0$ следующим образом.
Пусть, $ M_0$ - полюс, $ Ox$ - произвольно направленная ось в
плоскости $\Gamma_{123}$, $Oy$ -перпендикулярная ей ось в
плоскости $\Gamma_{123}$, $Oz$ - ось, перпендикулярная плоскости
$\Gamma_{123}$. В этих координатах оператор V записывается
следующим образом:
$$V \textbf{x} = (\lambda_1(x \cos\phi - y sin\phi), \lambda_1(x \sin\phi + y \cos\phi),\lambda_2z) + O(x^2 + y^2 + z^2)$$
 а обратное преобразование - в виде
 $$V^{-1}\textbf{x} = (\lambda^{-1}_1(x \cos\phi - y \sin\phi), \lambda^{-1}_1(x \sin\phi + y \cos\phi),\lambda^{-1}_2 z) +
 O(x^2 + y^2 + z^2).$$
Здесь $|\lambda_1|>1$ и $|\lambda_2|<1$. Зафиксируем некоторое
достаточно малое $\varepsilon.$ Без ограничения общности можем
считать что $(V^n \textbf{x})_4 \leq \varepsilon $ для
произвольного положи\-тельного  $n$. Рассмотрим окрестность $U$
точки $M_0 $ определенную в канонических координатах следующим
образом:
$$ U=\{x^2+y^2 \leq \delta ^2,0\leq z \leq \delta \}. $$
Константа $\delta $ выбирается таким образом, что выполняется
условие
$$ V^{-1}U \subset \{x^2 + y^2 \leq \delta^2 \}= U_1,$$
и при этом выбранная константа обеспечивает условие
$$ x_3= min\{x_1,x_3,x_4\}$$  для любого   $ \textbf{x}\in V^{-1}U.$
 После выбора $\delta $ ,выберем малое $\varepsilon$ так,что
 $$ U_1\bigcap \{x_3 \leq \varepsilon \} \subset U.$$
Теперь покажем,что для любого $n$, $ V^n \textbf{x} \in U.$
Рассмотрим последовательность $\{n_k\}$  такую, что $
V^{n_k}\textbf{x}\in U.$ Тогда
$$V^{n_k -1}\textbf{\textbf{x}}=V^{-1}(V^{n_k}\textbf{x}) \in V^{-1}U \subset U_1 ,$$
и так как
$$ V^{n_k -1}\textbf{\textbf{\textbf{x}}} \in \{x_4 \leq \varepsilon\}, $$
отсюда следует ,что $ V^{n_k -1}\textbf{x}\in U.$ В силу того,что
$n_k \rightarrow \infty$ при $ k\rightarrow \infty $,для любого
$k$ и для любого $t\leq n_k , V^t \textbf{x}\in U.$ Таким образом
$ \omega (\textbf{x})\subset U $ и так как $\omega
(\textbf{\textbf{x}})\subset \Gamma _{123},$ отсюда следует,что
существует достаточно малая окрестность $ W$(ее диаметр порядка
$const \cdot\delta$) точки $M_0$ в $\Gamma_{123}$ такая что
$\omega(\textbf{x})\subset W.$  Пусть $ \textbf{y} \in
\omega(\textbf{x})\neq M_0 $,тогда $V^n \textbf{y}
\in\omega(\textbf{\textbf{x}})$ и так как
$\omega(\textbf{\textbf{\textbf{x}}})$ замкнуто,то
$\omega(\textbf{y})\subset\omega(\textbf{\textbf{x}}).$ В то же
время $\omega(\textbf{y})\subset\partial\Gamma_{123}$ и
$V\bigcap\partial\Gamma_{123}=\emptyset.$ Из этого противоречия
следует,что $\omega(\textbf{\textbf{x}})=M_0.$ Применяя теорему
Гробмана-Хартмана [10] получим,что точка $\bf x$ лежит на
инвариантной кривой оператора $V,$ которая начинается в вершине
$(0,0,0,1)$ и заканчивается в точке $M_0.$ Этот случай не
представляет интереса.

Второй случай: $\omega$-предельное множество точки M не содержит
точку $ M_0$.Пусть $\textbf{y}\in\omega(\textbf{x})$ , тогда $V^{-n}
\textbf{y}\in\omega(\textbf{x})$ и так как $\omega(\textbf{x})$
замкнуто, то $\omega(\textbf{y})\subset\omega(\textbf{x}).$Легко
видеть что $\omega(\textbf{y})=M_0$ для всех
$\textbf{y}\notin\partial\Gamma_{123}.$ Поэтому, если такое
$\textbf{y}$ принадлежит $\omega({\textbf{x}}),$  то приходим к
противоречию с предположением, что $M_0$ не принадлежит
$\omega({\textbf{x}})$. Таким образом все $\omega $-предельные точки
$\textbf{x}$ принадлежат $\partial\Gamma_{123}.$ Итак, исключая
случай когда вектор $\textbf{x}$ принадлежит инвариантной кривой,
доказывается что $\omega(\textbf{x})\subset
\partial\Gamma_{123}.$
Это утверждение является точным аналогом соответствующего результата
для двумерного ксо [3].
\subsection{Матрицы третьего типа}
Запишем матрицу $A$  в канонической форме (5), где все числа $a_{ij}$
неотрицательны.

В этом случае отсутствуют притягивающие или отталкивающие вершины и
все вершины являются седловыми точками. Также внутри ребер и граней
$\Gamma_{234}$ и $\Gamma_{123}$ отсутствуют неподвижные точки ,а
внутри граней $\Gamma_{124}$ и $\Gamma_{134}$ найдутся по одной
неподвижной точке,причем одна из них седловая, а
другая-отталкивающая.

Их тип зависит от знака величины
$$I=-a_{12}a_{34}+a_{13}a_{24}+a_{14}a_{23}$$
{\bf Теорема 2.} Пусть  $\textbf{x}$- точка общего положения
,тогда
  $$
\omega(\textbf{x})\subset \Gamma_
{123}\cup\Gamma_{234}\cup\Gamma_{14}.$$ Исключением являются точки
границы симплекса и точки лежащие на инвариантной
кривой,начинающейся в отталкивающей неподвижной точке и
заканчивающейся в седловой точке, отличной от вершин.

Доказательство.При доказательстве существенно используется следующий
результат, доказанный в [6,8].

{\bf Предложение 1.}Пусть $\bf p$-отталкивающая неподвижная
точка.Тогда существует функция Ляпунова вида
$$ F(\textbf{x})=\prod_{i=1}^n x_i ^{p_i}.$$

Без ограничения общности можем считать, что неподвижная точка внутри
грани $\Gamma_{134}$ отталкивающая, а внутри грани $\Gamma_{124}$-
седловая.Тогда отталкивающая неподвижная точка имеет вид
$(p_1,0,p_3,p_4),$ а соответствующая функция Ляпунова имеет вид
$$ F(\textbf{x})=x_1^{p_1}x_3^{p_3}x_4^{p_4}.$$
Легко видеть, что вдоль траектории эта функция стремится к нулю.
Утверждение Предложения 1 можно переписать в следующей форме
$$ \lim_{n\rightarrow\infty} min\{(V^n \textbf{x})_1,(V^n \textbf{x})_3,(V^n \textbf{x})_4 \}=0
.$$ Эта форма более удобна для приложений.

Рассмотрим следующие случаи.

Первый случай. Пусть $\omega(\textbf{x})$ содержит седловую точку
$M_0$ внутри грани $\Gamma _{124}.$ Тогда для любой окрестности $U$
точки $M_0$ существует последовательность $\{n_k\}$ стремящаяся к
бесконечности такая, что $V^{n_k}\textbf{x}\in U.$ Введем
канонические координаты в окрестности $ U$ следующим образом. Пусть,
$ M_0$- полюс, $ Ox$ - произвольно направленная ось в плоскости
$\Gamma_{124}$, $Oy$ -перпендикулярная ей ось в плоскости
$\Gamma_{124}$, $Oz$ - ось, перпендикулярная плоскости
$\Gamma_{124}$. Тогда оператор V записывается следующим образом:
$$V \textbf{x} = (\lambda_1(x \cos\phi - y sin\phi), \lambda_1(x \sin\phi +
y \cos\phi),\lambda_2z) + O(x^2 + y^2 + z^2)$$
 а обратное преобразование - в виде
 $$V^{-1}\textbf{x} = (\lambda^{-1}_1(x \cos\phi - y \sin\phi),
 \lambda^{-1}_1(x \sin\phi + y \cos\phi),\lambda^{-1}_2 z) +
 O(x^2 + y^2 + z^2).$$
Здесь $|\lambda_1|>1$ и $|\lambda_2|<1$. Зафиксируем некоторое
достаточно малое $\varepsilon.$ Без ограничения общности можем
считать, что $(V^n \textbf{\textbf{x}})_4 \leq \varepsilon $ для
произвольного положи\-тельного  $n$. Рассмотрим окрестность $U$
точки $M_0 $ определенную в канонических координатах следующим
образом:
$$ U=\{x^2+y^2 \leq \delta ^2,0\leq z \leq \delta.\} $$
Константа $\delta $ выбирается таким образом что выполняется условие
$$ V^{-1}U \subset \{x^2 + y^2 \leq \delta^2 \}= U_1,$$
и при этом выбранная константа обеспечивает условие
$$ x_3= min\{x_1,x_3,x_4\}$$  для любого   $ \textbf{x}\in V^{-1}U.$
 После выбора $\delta $ ,выберем малое $\varepsilon$ так,что
 $$ U_1\bigcap \{x_3 \leq \varepsilon \} \subset U.$$
Теперь покажем,что для любого $n,$ $ V^n \textbf{x} \in U.$
Рассмотрим последовательность $\{n_k\}$  такую, что $
V^{n_k}\textbf{x}\in U.$ Тогда
$$V^{n_k -1}\textbf{x}=V^{-1}(V^{n_k}\textbf{x}) \in V^{-1}U \subset U_1 $$
и так как
$$ V^{n_k -1}\textbf{x} \in \{x_3 \leq \varepsilon\}, $$
отсюда следует ,что $ V^{n_k -1}\textbf{x}\in U.$ В силу того,что
$n_k \rightarrow \infty$ при $ k\rightarrow \infty $,для любого
$k$ и для любого $t\leq n_k , V^t \textbf{x}\in U.$ Таким
образом $ \omega (\textbf{x})\subset U $ и так как $\omega
(\textbf{x})\subset
\partial S^{4-1}$ отсюда следует,что существует достаточно малая
окрестность $ W$ (ее диаметр порядка $const \cdot\delta$) точки
$M_0$ в  $\Gamma_{123}$ такая что
$\omega(\textbf{\textbf{x}})\subset W.$ Пусть $ \textbf{y} \in
\omega(\textbf{x})\neq M_0 $,тогда  $V^n \textbf{\textbf{y}}
\in\omega(\textbf{x})$ и так как $\omega(\textbf{x})$ замкнуто,то
$\omega(\textbf{y})\subset\omega(\textbf{x)}.$ В то же время
$\omega(\textbf{y})\subset\partial\Gamma_{124}$ и
$V\bigcap\partial\Gamma_{123}=\emptyset.$ Из этого противоречия
следует,что $\omega(\textbf{x})=M_0.$ Применяя теорему
Гробмана-Хартмана [10] получим,что вектор $\bf x$ лежит на
инвариантной кривой оператора $V,$ которая начинается в вершине
$(p_1,0,p_3,p_4)$ и оканчивается в точке $M_0.$ Этот случай не
представляет интереса.

Второй случай:$\omega$-предельное множество точки $\textbf{x}$ не
содержит точки $ M_0$.Пусть
$$\textbf{y}\in\omega(\textbf{x})\cap \Gamma_{124}.$$
Тогда $V^{-n} \textbf{y}\in\omega(\textbf{x})$ и так как
$\omega(\textbf{x})$ замкнуто, то
$\omega(\textbf{y})\subset\omega(\textbf{x}).$Легко видеть что
$\omega(\textbf{y})=M_0$ для всех
$\textbf{y}\notin\partial\Gamma_{124}.$Поэтому, если такое
$\textbf{y}$ принадлежит $\omega({\textbf{x}})$ , то приходим к
противоречию с предположением, что $M_0$ не принадлежит
$\omega({\textbf{x}})$. Таким образом
$$\omega(\textbf{x})\subset \{x_1x_2x_4=0\}.$$
Аналогично
$$\omega(\textbf{x})\subset \{x_1x_3x_4=0\},$$
откуда
$$\omega(\textbf{x})\subset \{x_1x_2x_4=0\}\cap\{x_1x_3x_4=0\}=
\Gamma_{123}\cup \Gamma_{234}\cup\Gamma_{14}.$$ Теорема доказана.
\section{Геометрические и эргодические свойства}
В случае матриц первого типа траектория любой точки сходится к одной
из вершин.Если матрица принадлежит второму классу, то проблема
редуцируется к двухмерному ксо изученному в [3].Ниже будем
рассматривать только случай матриц принадлежащих третьему классу.
Для простоты будем предполагать, что числа $a_{ij}$ отличны от $0$ и
$1.$
\subsection{Геометрическая теорема}
Теорему 2, доказанную выше,можно переформулировать следующим
образом: если $ \varphi(\textbf{x})=max \{x_1 x_2x_4,x_1x_3x_4\}, $
тогда $\varphi(V^n \textbf{x})\rightarrow 0.$ Основным результатом
этой работы является следующая теорема.

{\bf Теорема 3.}Для матриц третьего класса соответствующий ксо
неэргодичен.

Предварительно докажем следующее предложение.

 {\bf Предложение 2.}Зафиксируем достаточно малое значение $\varepsilon$ и
 определим множество
 $$ U_{\varepsilon}=\{\textbf{x}\in S^{4-1}:x_2,x_3,x_4\leq\varepsilon\} .$$
Существуют константы $A$ и $B$  такие, что если $\textbf{x}\notin
U_{\varepsilon},V\textbf{x},\cdots ,V^{n-1}\textbf{x} \in
U_{\varepsilon},$ $V^n \textbf{x}\notin U_{\varepsilon},$ то
$$ n\geq A \log(\frac{B}{\varphi (\textbf{x})}).$$
 Доказательство.Легко видеть, что
 $$  \frac{(V^n \textbf{x})_4}{x_4}=\prod_{i=0}^{n-1} \frac{(V^{i+1}\textbf{x})_4}{(V^i
 \textbf{x})_4}\geq 2^n .$$
 Пусть $ \textbf{y}=V^{n-1}\textbf{x} \in U_{\varepsilon}.$ Тогда
 \begin{eqnarray*} (V^n \textbf{x})_2&=&(V\textbf{y})_2=y_2 (1-a_{12}y_1+a_{23}y_3+a_{24}y_4)\leq
 y_1(1-a_{12}(1-3\varepsilon)+2\varepsilon) \\
& \leq &y_2(1-a_{12}+5\varepsilon)\leq y_2\leq \varepsilon
\end{eqnarray*}
и
\begin{eqnarray*} (V^n \textbf{x})_3&=&(Vy)_3=y_3 (1-a_{13}y_1-a_{23}y_2+a_{34}y_4)\leq
 y_3(1-a_{13}y_1+\varepsilon)\\
&\leq& y_3(1-a_{13}(1-3\varepsilon)+\varepsilon)\leq
 y_3(1-a_{13}+4\varepsilon) \leq  y_3 \leq \varepsilon
\end{eqnarray*}
 Так как $ V^n \textbf{x}\notin U_{\varepsilon} $ ,то $(V^n \textbf{x})_4\geq \varepsilon $
 откуда следует
$$ 2^n \geq  \frac{\varepsilon }{x_4}= \varepsilon
\frac{x_1 max\{x_2,x_3\}}{\varphi (\textbf{x})}.$$

В силу следующего неравенства
$$ (V\textbf{x})_2= x_3(1-a_{12}x_1+a_{23}x_3+ a_{24}x_4)\geq
x_2(1-a_{12})$$ имеем
$$x_2\leq \frac{\varepsilon}{1-a_{12}}. \eqno (6)$$
Аналогично
$$(V\textbf{x})_3=x_3(1-a_{13}x_1-a_{23}x_2+a_{34} x_4)\geq x_3
(1-max\{a_{13},a_{23}\}),$$ так что
$$x_3 \leq \frac{\varepsilon}{(1-max\{a_{13},a_{23}\})}, \eqno (7)$$
и
$$(V\textbf{x})_4=x_4(1+a_{14}x_1-a_{24}x_2 -a_{34} x_4)\geq x_4
(1-max\{a_{24},a_{34}\}),$$ откуда
$$ x_4\leq \frac{\varepsilon}{1-max\{a_{24},a_{34}\}}. \eqno (8)$$
Все эти неравенства (6-8)можем записать как
$$ x_2,x_3,x_4 \leq C_{abs} \varepsilon ,$$
где $C_{abs}$-абсолютная константа зависящая только от $ a_{ij}.$
Тогда $$ x_1 \geq 1-3C_{abs}\varepsilon.$$ Теперь

\begin{eqnarray*}(V\textbf{x})_4&=&x_4(1+a_{14}x_1-a_{24}x_2-a_{34}x_3)\\
&\geq&
x_4(1+a_{14}(1-3C_{abs}\varepsilon)-(a_{24}+a_{34})C_{abs}\varepsilon)\\
&\geq& x_4(1+a_{14}
-C_{abs}(3a_{14} +a_{24}+a_{34})\varepsilon)\geq x_4,
\end{eqnarray*} откуда
$x_4\leq \varepsilon $ и $ max\{x_2,x_3\}\geq\varepsilon. $ Таким
образом
$$ 2^n\geq \frac{\varepsilon^2 (1-3C_{abs}\varepsilon)}{\varphi(\textbf{x})}.$$
Предложение доказано.

В качестве следствия имеем следующие теоремы

{\bf Теорема 4.}Время, проведенное точкой $\textbf{x}$ в окрестности
вершины $(1,0,0,0)$ за один сеанс, стремится к бесконечности.

{\bf Теорема 5.}Среднее время, проведенное точкой $\textbf{x}$ вне
$\varepsilon$-окрестности всех вершин, равно $0$.

Доказательство. Ясно,что точка движется по одному из следующих
маршрутов:
$$ 1\rightarrow 4\rightarrow 3\rightarrow 2\rightarrow1$$
или
$$ 1\rightarrow 4\rightarrow 2\rightarrow1$$
или
$$ 1\rightarrow 4\rightarrow 3\rightarrow1.$$
Здесь $1,2,3,4$ означает   $\varepsilon$-окрестность
1-ой(соответственно 2-ой,3-ей,4-ой ) верши\-ны .Очевидно,время
необходимое для перемещения из окрестности одной вершины в
окрестность другой, равномерно ограничено.Обозначим верхнюю грань
таких времен через $T_{abs}.$ Каждый цикл движения включает в себя
один сеанс посещения первой окрестности и не более 4-х сеансов
перемещения из одной окрестности в другую.

Зафиксируем прозвольно большое $K$.В силу предыдущей теоремы
,сеансы посещения первой окрестности длятся как минимум время
$K$$К,$
 то есть асимптотическое число циклов за время $n$ не превысит
$\frac{n}{K}$.Время потраченное на перемещение из одной
окрестности в другую не превысит $\frac{4T_{abs} n}{K}$.Таким
образом ,среднее время пребывания вне окрестностей вершин не более
,чем  $\frac{4T_{abs} }{K}$.

Так как  $K$ произвольно,отсюда следует утверждение теоремы.
\section{Построение функций Ляпунова}
В заключение укажем способ построения функций Ляпунова в дополнение
к уже известным.Новые функции Ляпунова имеют вид
$$ F(\textbf{x})=\prod_{i=1}^4 x_i^{\lambda_i}.$$
Заметим,что здесь не требуется положительности чисел $\lambda_i$.
Очевидно
$$ F(V\textbf{x})=F(\textbf{x})\prod_{i=1}^4 (1+(A\textbf{x})_i)^{\lambda_i}=
F(\textbf{x})G(\textbf{x}).$$

{\bf Предложение 3.}Если $G<1$ во всех вершинах,тогда $F$-функция
Ляпунова.

Доказательство.Если $G<1$ во всех вершинах, то существует
окрестность $U$ вершин такая ,что $G<q<1$ в $U$. Обозначим
максимальное значение $G$ через $L$. Пусть $\tau_n$ -время
проведенное вне $U$ за общее время $n$. Тогда
 $$F(V^n \textbf{x})=F(\textbf{x})\prod _{i=0}^{n-1} G(V^i \textbf{x})
\leq F(\textbf{x}) q^{n-\tau_n} L^{\tau_n}\rightarrow 0,$$
 так как среднее время проведенное вне
 $U$ ,то есть  $\frac{\tau_n}{n},$ стремится к нулю.
 Рассмотрим условие предложения.Ее можно переписать в виде
   $$ \prod_{i=1}^4 (1+\tilde{a}_{ij})^{\lambda_i}<1.$$
Логарифмируя, имеем
   $$\sum_{i=1}^4 \log(1-\tilde{a}_{ji})\lambda_i <0.$$
Если определить матрицу
   $$ B=(b_{ij})=(-\log(1-\tilde{a}_{ij})),$$
 то это условие запишется как
   $$B {\bf \lambda } >0$$
   Это значит, что разбиение пространства параметров связано  с некоторыми характеристиками матрицы
   $B$.

{\bf Литература}

1. Бернштейн С.Н.,  Решение одной математической проблемы,
связанной с теорией наследованности, {\it Уч. Зап. Научно-Исслед.
каф. Укр. отд. Мат.,} {\bf 1} : 83-115 (1924).

2. Ганиходжаев Н.Н., Применение теории гиббсовских мер к
математической генетике, {\it ДАН Россия.} {\bf 61}: 321-323
(2000).

3. Ганиходжаев Н.Н.,Занин Д.В.,  Об одном необходимом условии эргодичности
квадратичных операторов определенных на двумерном симплексе, {\it Успехи

Математических наук } {\bf 59} № 3: 571-572 (2004).

4. Ganikhodjaev N.N., Rozikov U.A.,  On quadratic stochastic
operators generated by Gibbs distributions, {\it Regular and
Chaotic Dynamics.} {\bf 11}: No. 3 (2006).

5. Ганиходжаев Р.Н., Семейство квадратичных стохастических
операторов действующих на $S^2$. {\it ДАН РУз.} № 1: 3-5 (1989).

6. Ганиходжаев Р.Н., Квадратичные стохастические операторы,
функция Ляпунова и турниры, {\it Мат. Сб.} {\bf т.83.} № 8:
121-140 (1992).

7. Ганиходжаев Р.Н., Об определении бистохастического квадратичного
оператора. {\it УМН.} {\bf 48}, № 4: 244-246 (1992).

8. Ганиходжаев Р.Н., Карта неподвижных точек и функции Ляпунова
для одного класса дискретных динамических систем. {\it Мат.
Заметки} {\bf 56}: 1125-1131 (1994).

9. Захаревич М.И., О поведении траекторий и эргодической гипотезе
для квадратичных отображений симплекса. {\it УМН}, {\bf 33}:
207-208 (1978).

10. Хартман П , Обыкновенные дифференциальные уравнения Москва,Мир
(1970)

11. Kesten H.,Quadratic transformations:a model for population
growth,I,II, {\it Adv. Appl.Prob.}{\bf }:1-82,179-228 (1970)

12. Lyubich Yu.I. Mathematical structures in population genetics.
{\it Biomathematics}, {\bf 22,} Springer-Verlag, 1992.

13. Stein P.R., Ulam S.M. Nonlinear transformations studies on
electronic computers. {\it Rozprawy Mat.} {\bf 39}: 1-15 (1964).

14. Ulam S.M., {\it A collection of mathematical problems } ,
Interscience Publishers, New-York-London (1960)

\end{document}